\UseRawInputEncoding
\documentclass[12pt, leqno]{article}
\usepackage{latexsym}
\usepackage{enumerate,amsmath,amsthm,amssymb,amsfonts,amscd}

\newtheorem{defn}{Definition}[section]
\newtheorem{thm}[defn]{Theorem}
\newtheorem{cor}[defn]{Corollary}

\numberwithin{equation}{section}

\def\O{\mbox{${\cal O}$}}

\newcommand{\cc}{\mathbb C}

\renewcommand{\a}{\alpha}

\input{arrow.tex}

\begin{document}
\title{\bf Corona and Gleason's problems for certain Banach spaces of holomorphic functions}

\author{S. R. PATEL (https://orcid.org/0000-0003-4960-3171)}

\maketitle \vspace{1.5cm}

\noindent Department of Mathematics, C. U. Shah University, Wadhwan City,\\
Gujarat, INDIA\\
E-mail: srpatel.math@gmail.com
\begin{tabbing}
\vspace{1.5cm} \noindent
{\bf 2020 Mathematics Subject Classification:} \= Primary 46E15;\\
\> Secondary 30H80, 32A10\\
\> 32A37, 32A70
\end{tabbing}
\noindent {\bf Abstract}-- We use a combination of the classical
dead-end space method and the interpolation (by bounded
holomorphic functions) method to affirmatively solve Corona and
Gleason's problems for familiar Banach spaces of holomorphic
functions such as (Besov-)Lipschitz spaces of holomorphic
functions on bounded symmetric domains (including the unit ball
and polydisc) in $\cc^n$, holomorphic mean (Besov-)Lipschitz
spaces, (Besov-)Hardy-Sobolev spaces and some special classes of
holomorphic mixed-norm (Bergman) spaces (including weighted
Bergman spaces) as well as some special classes of holomorphic
(mean) Besov-Lipschitz spaces on the unit ball in $\cc^n$. Indeed,
we extend our recent work on both problems for holomorphic Lipschitz algebras.\\

\noindent {\bf Keywords:} Corona problem, Gleason's problem,
Banach spaces of holomorphic functions, dead-end space and
interpolation methods, holomorphic Lipschitz algebras
\newpage
\section{Introduction} The objective of this short
article is to provide affirmative solutions to the most
prestigious Corona problem (1962) and the well-known Gleason's
problem (1964) for certain familiar Banach spaces of holomorphic
functions on specific domains in $\cc^n$. We first introduce both
problems with a brief history suitable to our work.

\noindent {\bf Corona problem.} In 1942, Kakutani asked the
following fundamental question: Let $H^\infty(D)$ be the Banach
algebra of all bounded, holomorphic functions on the open unit
disc $D$ in the complex plane $\cc$, equipped with the usual
sup-norm, and let $f_j \,\in\,H^\infty(D), \,j = 1, 2,\dots, n$,
with the property that $\sum_{j = 1}^n|f_j(z)|\,>\,\epsilon$ for
some $\epsilon\,>\,0$ and all $z\,\in\,D$ (namely, a finite set of
corona data on the disc $D$), do there exist $g_j
\,\in\,H^\infty(D), \,j = 1, 2,\dots, n$, such that $\sum_{j =
1}^nf_j(z)g_j(z)\,=\,1$ for all $z\,\in\,D$? This latter equation
is known as the B\'{e}zout equation ~\cite[Theorem V.1.8]{14}.

The equivalent formulation of this question is to ask whether $D$
is dense in the maximal ideal space (aka spectrum)
$M(H^\infty(D))$ w.r.t. the relative weak*-topology (equivalently,
whether the corona is empty). Carleson affirmatively solved the
problem on the disc $D$ ~\cite{4}.

By a {\it domain} we mean an open, connected and bounded set $G$
in $\cc^n$. At this time there is no known domain in $\cc$ on
which this problem is failed.

Since 1962, the corona problem is open for arbitrary domains
(including the balls and polydiscs) in $\cc^n,\,n\,>\,1$, and is a
lot more complicated. There have been investigations of this
problem on various types of domains, as well as in various Banach
spaces of holomorphic functions on these domains. For example, see
some partial results on this problem in ~\cite{3, 28, 35}. Not
only this, but a variety of counterexamples to the corona problem
in the theory of SCV have been produced in ~\cite{13, 33}. The
extensive exposition of history of the research of this problem is
outlined in ~\cite{10}.

In fact, after the Carleson's solution, most of the subsequent
attempts to extend it for various Banach spaces of holomorphic
functions on more general domains in $\cc^n$ were aimed at solving
the B\'{e}zout's equation by various classical methods (mostly, by
solving the $\bar{\partial}-$equation); see, e.g., ~\cite{35}.
These approaches required hard classical analytical methods and
has met serious difficulties in the case when $n > 1$, {\bf not}
leading to the solution to the original $H^\infty$ corona problem.
For example, Varopoulos in ~\cite{35} realized that his result
{\it does not imply} that the corona problem fails in several
variables-only that the $\bar{\partial}-$technique with that {\it
particular definition of Carleson measure} fails. We discuss the
functional analytic methods different from those classical
analytic methods. In particular, we have established an
equivalency between the corona problem and Gleason's problem
(which we describe briefly as follows) for Lipschitz algebras of
holomorphic functions in Theorems 2.2 and 2.3 below (see
~\cite{28} for details).

\noindent {\bf Gleason's problem.} Let $A(G)$ be the Banach
algebra of all holomorphic functions on a domain $G$ in $\cc^n$
that extend by continuity to the closure of $G$, equipped with the
sup-norm. Thus, it is a closed subalgebra of $H^\infty(G)$. The
Gleason's problem is to decide whether the ideal in $A(G)$ (or
$H^\infty(G)$), consisting of functions vanishing at a fix point
$\lambda = (\lambda_1, \lambda_2, \dots, \lambda_n) \in G$, is
algebraically finitely generated by the shifted coordinate
functions $z_j - \lambda_j$ for $j = 1, 2, \dots, n$ (see
~\cite{1} and other references therein for details). The
equivalent algebraic formulation is: suppose $f\, \in\, A(G)$ (or
$H^\infty(G)$) is such that $f(\lambda)\, =\, 0$ for some
$\lambda\, \in\, G$. Then there exists $f_1,\,f_2,\,\dots,\,f_n$
in $A(G)$ (or $H^\infty(G)$) so that $f\,=\,\sum_{j = 1}^n(z_j -
\lambda_j)f_j$. This formulation allows us to consider ``certain''
Banach spaces of holomorphic functions without having Banach
algebra structures; see Theorem 3.1 below. This is actually a
variant of the original Gleason's problem, posed for $A(G)$, where
$G$ is the unit ball in $\cc^2$ and $\lambda = (0, 0)$ ~\cite{15}.
The positive solution for strictly pseudoconvex domains with
smooth boundaries was established in ~\cite{18}. It has been even
strengthened in ~\cite{20} as follows: The continuous linear
operators
$$T_j : A(G) \rightarrow
A(G),\;f\,\mapsto\,T_j(f),\;j\,=\,1,\,2,\,\dots,\,n$$ solving
Gleason's problem exist, which means that $$f(z) -
f(\lambda)\,=\,\sum_{j = 1}^n(z_j -
\lambda_j)(T_j(f))(z),\;z\,\in\,G, f\,\in\,A(G).$$

This implies, in particular, that the part of the spectrum
$M(A(G))$ in the fibers over $G$ is mapped in a bijective way to
$G$ by the Gel'fand transform $(\hat{z}^1,\,\dots,\,\hat{z}^n)$ of
the coordinate functions. By the result of ~\cite{16}, the entire
spectrum $M(A(G))$ is mapped bijectively onto $\overline{G}$ where
$G$ is a smooth pseudoconvex domain in $\cc^n$ (it also suffices
to consider a pseudoconvex domain $G$ with a Stein neighbourhood
basis; see ~\cite{31}). Hence the Gel'fand topology coincides with
the Euclidean topology on $\overline{G}$ and $G$ is open in
$M(A(G))$.

A domain $G$ has the Gleason $A-$property if Gleason's problem has
an affirmative solution at all points of $G$ for a locally Stein
algebra $A$ ~\cite{27}. Clos calls it the Gleason solvability
property of a domain $G$ for $H^\infty(G)$ ~\cite{9}. The Gleason
solvability was first proved by Leibenson, albeit informally
~\cite{32}. The Gleason solvability for holomorphic Bergman
spaces, holomorphic (harmonic) mixed-norm spaces, Bergman-Sobolev
spaces and the holomorphic (harmonic) Bloch spaces were studied in
~\cite{32, 24, 26, 37} (~\cite{19}). The domains considered were
initially the unit ball, but then were generalized to (bounded)
strongly pseudo(convex) domains with $C^2$ boundary. There is also
considerable interest on other function spaces and domains,
including the weighted Bergman spaces on egg shaped domains, as
seen in ~\cite{30}. The Gleason solvability is sometimes referred
to as Hefer's condition (as developed by Hefer in 1940). Many of
the methods used in the study of Hefer's condition involve the
integral representation techniques developed in ~\cite{17, 29}.
\section{Past work as a background}
This paper contains a continuation of the work begun in
~\cite{28}, and we shall feel free to use the terminology and
conventions established there. In ~\cite{27}, we extend the
Gleason's result to finitely generated ideals in Fr\'{e}chet
algebras, answering affirmatively a question posed in ~\cite{5}
for Fr\'{e}chet algebras. As consequences, locally Stein algebras
are completely characterized by intrinsic properties within the
class of Fr\'{e}chet algebras, using which an affirmative answer
to Gleason's problem for such algebras is provided. In fact, we
remark that complex analysts may find this functional analytic
approach interesting from applications point of view. For example,
see ~\cite[Theorems 4.2, 4.3 and Corollary 4.4]{27}, and
references (and their reviews) to ~\cite{15} in MathSciNet
(MR0159241 (28 \#2458)). As we will see, this is, indeed, the
case. Banach algebras satisfying ~\cite[Lemma 3.1]{27} (in
particular, the polydisc algebra $A(D^n)$, the ball algebra
$A(B^n)$ and $H^\infty(G)$) are locally Stein-Banach algebras. In
~\cite{28}, we show that Banach subalgebras $A^k(\overline{G})$
and $\textrm{Lip}_H^{k,\,\a}(G)$ of $H^\infty(G)$, whose Gleason
solvability we have established, are also locally Stein-Banach
algebras.

In ~\cite[Corollary 4.4]{27}, we show that the dense open subspace
has the Gleason $A$-property, where $A$ is a semi-simple locally
Stein algebra; in particular, $A$ can be either $H^\infty(G)$ or
$A(G)$, where $G$ is a (strictly or weakly) pseudoconvex domain in
$\cc^n$ with various boundary conditions and containing the
origin. Thus the abstract functional analytic method recaptures
the classical results on Gleason's problem obtained by Beatrous
Jr. ~\cite{2}, Forn\ae ss and \O vrelid ~\cite{12}, Kerzman and
Nagel ~\cite{22}, Lieb ~\cite{23}, Noell ~\cite{25}, and Backlund
and F\"{a}llstr\"{o}m ~\cite{1} (see other references therein for
papers on this problem). This problem has recently been solved for
Lipschitz algebras of holomorphic functions. In fact, we
affirmatively answer a question of Forn\ae ss and \O vrelid (1983)
for holomorphic H\"{o}lder and Lipschitz spaces (~\cite[Remarks E.
2]{27}).

We found it necessary (mathematically, ``sufficient") to consider
the Gleason solvability while solving the corona problem. On the
other hand, although we did not spell it out explicitly, but
considered the corona theorem as a hypothesis while solving
Gleason's problem ~\cite[Theorem 4.3 and Corollary 4.4]{27}. In
fact, \.Zelazko established a pleasant result ~\cite[Proposition
2]{36}, showing a possible connection between these two problems
through their algebraic formulations. Moreover Kerzman and Nagel
used the Koszul complex while solving Gleason's problem for
certain function algebras on certain domains in $\cc^n$. All this
developments stemmed to think on an equivalency of these two
problems, as well as to gave a short elegant proof, by using the
(maximal) ideal structures of Banach algebras $A^k(\overline{G})$
and $\textrm{Lip}_H^{k,\,\a}(G,\,d)$ alone in ~\cite{28}. Thus we
kept ourselves away from the analytical aspects such as Gleason
parts of the spectra of these algebras.

In this paper, using this equivalency as the main source, we apply
the {\bf dead-end space} method and the {\bf interpolation} method
in next sections, in order to obtain positive results for both
problems for familiar Banach spaces of holomorphic functions on
bounded symmetric domains in $\cc^n$ (including the unit ball and
polydisc). Thus the solvability of both problems for Lipschitz
algebras of holomorphic functions plays a crucial role.

Let $G$ be a smooth strictly pseudoconvex domain in $\cc^n$
~\cite{34} or, else, $G$ be a pseudoconvex domain with a Stein
neighbourhood basis ~\cite{31}, e.g., Stein domains in $\cc^n$, in
the below results. First, we have the following crucial theorem,
extending results of ~\cite{9, 22}.
\begin{thm} \label{Thm 4_A(G)} Let $G$ be a ``pseudoconvex" domain
in $\cc^n$, and let $0\,\leq\,k\,<\,\infty$ and
$\a\,\in\,(0,\,1]$. Suppose $f\, \in\, \textrm{Lip}_H^{k,\,\a}(G,
\,d)$ (resp., $f\,\in\,A^k(\overline{G})$) and $f(\lambda)\, =\,
0$ for some $\lambda\, \in\, G$. Then there exists
$f_1,\,f_2,\,\dots,\,f_n$ in $\textrm{Lip}_H^{k,\,\a}(G, \,d)$
(resp., in $A^k(\overline{G})$) so that $f\,=\,\sum_{j = 1}^n(z_j
- \lambda_j)f_j$. In particular, $G$ is Gleason solvable for
$\textrm{Lip}_H^{k,\,\a}(G, \,d)$ and $A^k(\overline{G})$.
\end{thm}
\noindent {\it Proof.} For the proof, follow ~\cite[Theorem
3.2]{28}. $\hfill \Box$\\

The above theorem holds true if $G$ is a ``polydomain" (including
the polydisc) in $\cc^n$, obtained as a product of smooth strictly
pseudoconvex domains (see ~\cite[Theorem 3.3]{28}). As a part of
our goal to establish an equivalency between the two problems, we
have the following result from the above theorems.
\begin{thm} \label{Corona} Let $G$ be a ``pseudoconvex" domain or else a ``polydomain" in $\cc^n$. Then the corona problem holds true in
$A^k(\overline{G})$ and $\textrm{Lip}_H^{k,\,\a}(G,\,d)$, where
$0\,\leq\,k\,<\,\infty$ and $\a\,\in\,(0,\,1]$. That is, the
spectrum $M(A^k(\overline{G}))$ (resp.,
$M(\textrm{Lip}_H^{k,\,\a}(G,\,d))$ for each
$0\,\leq\,k\,<\,\infty$ and for each $\a\,\in\,(0,\,1]$) is the
Gel'fand closure of (the canonical image of) $G$.
\end{thm}
\noindent {\it Proof.} For the proof, follow ~\cite[Theorem
3.4]{28}. $\hfill \Box$\\

Then, in the converse direction, we have the following
\begin{thm} \label{Gleason} Let $G$ be a
relatively compact ``pseudoconvex" domain or else a ``polydomain"
in $\cc^n$. Then $G$ is Gleason solvable for $A^k(\overline{G})$
and $\textrm{Lip}_H^{k,\,\a}(G,\,d)$ for $0\,\leq\,k\,<\,\infty$
and $\a\,\in\,(0,\,1]$, provided that the corona problem holds
true in these algebras.
\end{thm}
\noindent {\it Proof.} For the proof, follow ~\cite[Theorem
3.6]{28}. $\hfill \Box$\\

As a corollary to Theorem 2.2, since the ball and polydisc are
(essentially the only) basic examples of a pseudoconvex/poly
domain in $\cc^n$ of our kind, we have the following
\begin{cor} \label{Corona_Special} Let $G$ be either an open unit ball or a polydisc in
$\cc^n$. Then the corona problem holds true in $A^k(\overline{G})$
and $\textrm{Lip}_H^{k, \a}(G, d)$, where $0 \leq k < \infty$ and
$\a \in (0, 1]$. $\hfill \Box$
\end{cor}
\section{Main Results}
We next use the Gleason solvability for holomorphic Lipschitz
algebras to further solve this problem for other familiar Banach
spaces of holomorphic functions such as Lipschitz spaces of
holomorphic functions over bounded symmetric domains in $\cc^n$
~\cite{6}, holomorphic mean Lipschitz spaces and Hardy-Sobolev
spaces on the unit ball in $\cc^n$ ~\cite{7}, and some special
classes of mixed-norm (Bergman) spaces (including weighted Bergman
spaces) of holomorphic functions on the unit ball in $\cc^n$
~\cite{21}. We remark that in the net of several Banach spaces of
holomorphic functions, analysts gave more than one name to one
such Banach space of holomorphic functions. For example, what Cho
and Zhu called special classes of mixed-norm spaces (e.g.,
holomorphic (mean) Lipschitz spaces) are exactly same as special
classes of the mixed-norm Bergman spaces (e.g., holomorphic (mean)
Besov-Lipschitz spaces, named by Jevti\'{c} and Pavlovi\'{c}
~\cite{21}), which is clear from various characterizations
obtained for these spaces in ~\cite{7, 8, 21}. We here define them
for the reader's convenience.

We first recall the definition of Lipschitz spaces over bounded
symmetric domains in $\cc^n$, taken from ~\cite[\S 1, 2]{6}, as
follows. Let $0\,<\,p,\, q\,\leq\,\infty$, and $\a$ be real. Then
$H\Lambda(\a, \,p,\, q)$ is the complex linear space of
holomorphic functions $f$ on a bounded symmetric domain $D$ such
that
$$\|f\|_{\a,\,p,\,q}\,=\,\mathcal{N}_{p,\,q,\,1}(J^{n(-\a-1)}f)\,<\,\infty,$$
where $J^\a$ is a multiplier transformation of $f$, and, for $\a
\,>\,0$, $$\mathcal{N}_{p,\,q,\,\a}(f)\,=\,(\int_0^1(1 - r)^{nq\a
- 1}\,M_p^q(f,\,r)\,dr)^{1/q}\;\;(0\,<\,q\,<\,\infty),$$
$$\mathcal{N}_{p,\,\infty,\,\a}(f)\,=\,\sup_{0\,\leq\,r\,<\,1}((1 -
r)^{n\a}\,M_p(f,\,r))$$ such that, for $0\,<\,r\,<\,1$, we set
$$M_p(f, \,r)\,=\,(\frac{1}{V}\int_{bdry\,D}|f(rt)|^p\,ds_t)^{1/p}\;\;(0\,<\,p\,<\,\infty),$$
$$M_{\infty}(f,\,r)\,=\,\sup_{t\,\in\,bdry\,D}|f(rt)|,\;\;\|f\|_p\,=\,\sup_{0\,\leq\,r\,<\,1}\,M_p(f,\,r)\,=\,\lim_{r\,\rightarrow\,1^-}\,M_p(f,\,r).$$
The class of such functions $f$ with $\|f\|_{p}\,<\,\infty$ is the
Hardy space $H^p(D)$. Now, we recall the definitions of
holomorphic mean Lipschitz spaces and Hardy-Sobolev spaces on the
unit ball $\mathbb{B}^n$ (canonical example of a bounded symmetric
domain) in $\cc^n$, taken from ~\cite{7}, as follows. For
$1\,\leq\,p\,<\,\infty$, $1\,\leq\,q\,\leq\,\infty$, and $\a
\,>\,0$, we consider the mixed-norm space $\Lambda_\a^{p,\,q}$
consisting of $f\,\in\,H(\mathbb{B}^n)$ such that $$(1 - r)^{k -
\a}\,M_p(\mathcal{R}^kf,\,r)\,\in\,L^q([0,\,1),\,\frac{dr}{1 -
r}),$$ where $\mathcal{R}^k$ is the special fractional radial
differential operator such that $k$ is any integer greater than
$\a$. When $q\,=\,\infty$, we will write $\Lambda_\a^p$ instead of
$\Lambda_\a^{p,\,\infty}$. When $p\,=\,q\,=\,\infty$, the space
$\Lambda_\a^{\infty,\,\infty}$ is exactly the classical
holomorphic Lipschitz space $\Lambda_\a\,=\,\textrm{Lip}_H^\a$, as
noted in ~\cite[\S 2]{28}. In fact, $J^\a$ and $\mathcal{R}^\a$
are inverse-like or equivalent operators on Lipschitz spaces over
bounded symmetric domains in $\cc^n$, providing equivalent norms
on these spaces (see ~\cite[Theorems D and 4]{6} and
~\cite[Theorem 3.12]{7}). In ~\cite[\S 1]{21}, they are called
Besov-Lipschitz spaces; also, when $1\,\leq\,p\,=\,q\,<\,\infty$,
we have
$$\|f\|_{p, p, \a}\,=\,\int_{\mathbb{B}^n}|f(z)|^p(1 - |z|^2)^{p\a
- 1}d\nu(z).$$ Therefore, $H^{p, p, \a}$ is a weighted Bergman
space. For $1\,\leq\,p\,\leq\,\infty$, and $\a$ real, the
Hardy-Sobolev space $H_\a^p$ consists of $f\,\in\,H(\mathbb{B}^n)$
such that $\mathcal{R}^{\a}f\,\in\,H^p$ ~\cite{7}. For
$f\,\in\,H_\a^p$, we write
$$\|f\|_{H_\a^p}\,=\,\|\mathcal{R}^{\a}f\|_{H^p}.$$ Then, for
$1 \leq q < \infty$, the Besov-Hardy-Sobolev space $H_\a^{p,\,q}$
consists of $f \in H(\mathbb{B}^n)$ such that $f \in
H_\a^{p}\,\bigcap\,\Lambda_\a^{p,\,q}$. For $f \in H_\a^{p,\,q}$,
by ~\cite[\S 1.1]{21}, we write
$$\|f\|_{H_\a^{p,\,q}}\,:=\,\max\{\|\mathcal{R}^{\a}f\|_{H_\a^{p}},\,\|\mathcal{R}^{s}f\|_{p,\,q,\,s
- \a}\}
\;\textrm{or}\;\|\mathcal{R}^{\a}f\|_{H_\a^{p}}\,+\,\|\mathcal{R}^{s}f\|_{p,\,q,\,s
- \a}.$$ We also consider some special classes of holomorphic
(mean) Besov-Lipschitz spaces such as the spaces $\Lambda_{n,
\phi}^p$ and the integrated spaces $\Lambda_{n, \phi}^{p, q}$ on
the unit ball in $\cc^n$, where $\phi$ is an increasing, positive,
continuous function on the interval $(0,\,1]$, from ~\cite[\S
1]{21} as follows. The space $\Lambda_{n, \phi}^p$ consists of
those $f\,\in\,H(\mathbb{B}^N)$ for which
$$M_p(r,\,\mathcal{R}^nf)\,=\,\O(\frac{\phi(1 - r)}{(1 -
r)^n}),\;r\,\rightarrow\,1^{-}.$$ The space $\textrm{Lip}_{n,
\phi}^p$ consists of those $f\,\in\,H(\mathbb{B}^N)$ for which
$$\omega_n^+(\delta,\,f)_p\,=\,\O(\phi(\delta)),\;\delta\,\rightarrow\,0^+.$$

Thus, we have the following two results in the below
\begin{thm} \label{Gleason_Solvability} Let $G$ be a bounded
domain in $\cc^n$. Then, for the following Banach spaces of
holomorphic functions on $G$: (A) $G$ is Gleason solvable, and (B)
the corona theorem holds true for $G$.
\begin{enumerate}
\item[{\rm (i)}] holomorphic (Besov-)Lipschitz spaces on a bounded
symmetric domain $G$ (including the unit ball and polydisc) in
$\cc^n$; \item[{\rm (ii)}] holomorphic mean (Besov-)Lipschitz
spaces and (Besov-)Hardy-Sobolev spaces on the unit ball in
$\cc^n$; \item[{\rm (iii)}] some special classes of holomorphic
mixed-norm (Bergman) spaces (including weighted Bergman spaces) on
the unit ball in $\cc^n$; and \item[{\rm (iv)}] some special
classes of holomorphic (mean) Besov-Lipschitz spaces such as the
spaces $\Lambda_{n, \phi}^p$ and the integrated spaces
$\Lambda_{n, \phi}^{p, q}$ on the unit ball in $\cc^n$, where
$\phi$ is an increasing, positive, continuous function on the
interval $(0,\,1]$.
\end{enumerate}
\end{thm}
\noindent {\bf Important remarks on main results.} We give some
comments on our functional analytic approach of a combination of
the dead-end space method and the interpolation (by bounded
holomorphic functions) to solve both problems for familiar Banach
spaces of holomorphic functions on bounded symmetric domains in
$\cc^n$. First, we note that due to the algebraic nature of both
problems, it is very much clear that one does not require the
Banach algebra structure of the underlying Banach space of
holomorphic functions for which we want to solve both problems;
equivalently, analysts did give solutions of both problems for
certain Banach spaces of holomorphic functions of specific domains
by applying various classical analytic methods [1]-[4], [9], [12],
[17]-[20], [22]-[26], [29], [30], [32], [35]. However, as we
remark in \S 1, those classical methods are hard and complicated
in the SCV case. So, in ~\cite{28}, we first establish both
problems for Lipschitz algebras $\textrm{Lip}_H^{k,\, \a}(G, d)$
of holomorphic function on certain domains in $\cc^n$ as well as
for Banach algebras $A^k(\overline{G})$ by applying a certain
functional analytic approach. This provides a concrete base for
the dead-end space method and the interpolation method to be used
in this paper, in order to discuss solutions of both problems for
familiar Banach spaces of holomorphic functions on bounded
symmetric domains in $\cc^n$. For more details on bounded
symmetric domains in $\cc^n$, see ~\cite[\S 1]{6}. In particular,
such domains are pseudoconvex domains with a Stein neighbourhood
basis as they are convex, circular and domain of holomorphy. In
fact, these domains are Reinhardt domains in $\cc^n$ or at least
have rotational symmetry. Further, the polydisc is a specific
example of a bounded symmetric domain.

Thus, our results are unique and very different in two ways: (1)
we move out from applying classical methods, to consider
functional analytic methods, which may have a broader scope in the
future; and (2) we do not repeat Banach spaces of holomorphic
functions for which solutions of both problems are available (in
other words, we do not recover the old results), as mentioned in
the above references.
\section{Proofs} Recall that Carleson in ~\cite{4} used the interpolation (by $H^\infty-$
functions) method to solve the corona problem for the open unit
disc in $\cc$. Dronov and Kaplitskii in ~\cite{11} used a
combination of the dead-end space method and the interpolation (by
bounded operators on the cones in certain Banach spaces) method to
solve the Pelczy\'{n}ski's basis problem from the theory of
nuclear Fr\'{e}chet spaces: whether each complemented subspace of
a nuclear Fr\'{e}chet space with basis also has a basis. Their
main result ~\cite[\S 3.3, Theorem 2]{11} is as follows: let $E$
be a nuclear K\"{o}the space with regular basis in the class
$(d_1)$, and let $F$ be a complemented subspace of $E$. Then $F$
has an absolute basis. The central idea of the dead-end space
method and the interpolation method is as follows: they construct
a system of elements to become a basis in $F$ as a common
orthogonal basis of two Hilbert spaces $F_0$ and $F_\infty$ such
that $F_\infty\,\subset\,F\,\subset\,F_0$ and $F_\infty$ is
embedded in $F_0$ compactly and densely (such a basis always
exits). Furthermore, taking appropriate weights $a_0(n)$ and
$a_\infty(n)$ they could always achieve
$F_0\,\subset\,E_0\,=\,\ell_2(a_0(n))$,
$F_\infty\,\subset\,E_\infty\,=\,\ell_2(a_\infty(n))$,
$E_\infty\,\subset\,E_r\,\subset\,E_0$, and for each r the triple
$(E_0, E_\infty, E_r)$ is an interpolation triple of Hilbert
spaces. The system of elements thus constructed is always a
complete minimal system in $F$, so it is a basis in $F$, provided
that the corresponding operators of partial sums $P_n$ are
equicontinuous. These operators are bounded on $P(E)\,=\,F$, where
$P$ is a continuous projection in $E$ that takes $E$ to $F$, but
they are unbounded in the two extreme spaces because $P$ is
continuous in the Fr\'{e}chet space $E$, but not continuous
(unbounded) in the extreme Hilbert spaces $E_0$ and $E_\infty$.
However, the operators $P_n$ are bounded on certain cones in $E_0$
and $E_\infty$, which are constructed from the operator $|P|$ in a
natural way. These cones have many nice properties (e.g., they are
lower semilattices), so they could use interpolation theorems for
operators bounded on them to show that the family $\{P_n\}$ is
equicontinuous.

We will use this combination, but for the interpolation (by
bounded holomorphic functions in certain Banach spaces) method to
establish the Gleason solvability and the corona theorem below.
Indeed, we use the dead-end space method to show the Gleason
solvability for these spaces of bounded symmetric domains $G$ in
$\cc^n$ (including the ball and the polydisc). The central idea of
this method is as follows.

We construct the triple $(F_0, F_\infty, F)$ of Banach spaces of
holomorphic functions of $G$ in $\cc^n$ so that $G$ is Gleason
solvable for $F_0$ and $F_\infty$ such that $F_\infty \subset F
\subset F_0$ continuously and $F_0 \cong F_\infty$ (compact and
dense embedding of $F_\infty$ in $F_0$ will also do). Furthermore,
taking appropriate Gleason solvable Banach spaces $E_0$ and
$E_\infty$ of holomorphic functions of $G$, we can always achieve
$F_0\,\subset\,E_0$ continuously, $F_\infty\,\subset\,E_\infty$
continuously, $E_\infty\,\subset\,E_r\,\subset\,E_0$ continuously
and for each $r$, the triple $(E_0,\, E_\infty,\, E_r)$ is an
interpolation triple of Gleason solvable Banach spaces of
holomorphic functions of $G$. The system of Banach spaces thus
constructed, always squeezes $F$. Thus, $F$ is Gleason solvable,
provided that the corresponding $F_0,\,F_\infty,\,E_0,\,E_\infty$
and $E_r$ are Gleason solvable. We remark that the dead-end space
method has been used to affirmatively solve certain well-known,
long-standing problems from (functional) analysis such as
Pelczy\'{n}ski's basis problem (Grothendieck's basis problem was
the most prestigious, general one among all basis problems in the
theory of nuclear Fr\'{e}chet spaces) ~\cite{11}.

Since we consider both the problems for bounded symmetric domains
in $\cc^n$ (including the unit ball and polydisc), we remark that
they are equivalent in view of Theorems 2.1, 2.2 and 2.4 above.
So, one can also prove Theorem 3.1 (B) first, and then, give proof
of Theorem 3.1 (A). \\

\noindent {\it Proof of Theorem 3.1 (A).} (i). We first remark
that for $\alpha\,\in\,(0, 1], \, k\,=\,0$ and $G$ being a bounded
symmetric domain (including the unit ball and polydisc) in
$\cc^n$, the Lipschitz algebras $\textrm{Lip}_H^\a(G,\,d)$,
considered in Theorem 2.1, are the classical holomorphic Lipschitz
spaces $\Lambda_\alpha\,=\,H\Lambda(\alpha/n,\,\infty,\,\infty)$,
discussed by Chen ~\cite[Theorem 10]{6} that are actually same as
the spaces $\Lambda_\alpha^{\infty,\,\infty}$, discussed by Cho
and Zhu ~\cite{7} for the unit ball. For $\a\,=\,1$, we here have
the special case of the Lipschitz algebra $\textrm{Lip}_H(G,\,d)$
of holomorphic functions, which is the space
$\Lambda_1\,=\,\Lambda_1^{\infty,\,\infty}$, discussed by Cho and
Zhu ~\cite{7}, that is actually same as the space
$H\Lambda(1/n,\,\infty,\,\infty)$, discussed by Chen
~\cite[Theorems 4 and 10]{6}). Moreover, by ~\cite[Lemma 2]{6}, we
see that a function $f \in H\Lambda(\alpha/n,\,\infty,\,\infty)$
is continuous on $G\,\bigcup\,b$, where $b$ is a Bergman-\u{S}ilov
boundary of $G$ (i.e., $f$ is uniformly continuous on $G$). For $1
< p < \infty$ and $q = \infty$, by ~\cite[Theorem 6]{6},
\begin{equation} \label{Thm 6} \Lambda_\alpha\,\subset\,H\Lambda(\alpha,\,p,\,\infty)\,\subset\,\Lambda_{\alpha - 1/p}
\end{equation}
(and these inclusions are continuous). Since $\Lambda_\alpha \cong
\Lambda_{\alpha - 1/p}$ by ~\cite[Theorem 4]{6}, we have that $G$
is Gleason solvable for $H\Lambda(\alpha,\,p,\,\infty)$ by the
dead-end space method and interpolation method, discussed above,
where $F_\infty = \Lambda_\a$ and $F_0 = \Lambda_{\a - 1/p}$ (note
that one suitably chooses $E_0 = \Lambda_\delta, \,E_\infty =
\Lambda_\beta$ and $E_r = \Lambda_\gamma$ to be the classical
holomorphic Lipschitz spaces (which are also Gleason solvable by
Theorem 2.1), because for $0 < \delta < \gamma < \beta < \a \leq
1$, we have $\Lambda_\a \subset \Lambda_\beta \subset
\Lambda_\gamma \subset \Lambda_\delta$ continuously). Along the
same lines, $G$ is Gleason solvable for
$H\Lambda(\alpha,\,\infty,\,q)$ for $1 \leq q < \infty$, since, by
~\cite[Theorem 5]{6}, for $\alpha \in (0,\,1),\;\beta = 1$, $p =
\infty$ and $1 \leq q < \infty$,
\begin{equation} \label{Thm 5}\Lambda_1\,\subset\,H\Lambda(\a,\,\infty,\,q)\,\subset\,\Lambda_\a
\end{equation}
(and these inclusions are continuous). The inclusion relation (2)
(also, (3) and (5) below) are the best examples of the spaces with
$q \,(\textrm{and}\,q_1) < \infty$ interpolating between the
spaces with $q \,(\textrm{and}\,q_2) = \infty$, hence a mention of
the interpolation method is appropriate. Since $J^\a$ and
$\mathcal{R}^\a$ provide equivalent norms on Lipschitz spaces, we
have that $G$ is Gleason solvable for holomorphic Besov-Lipschitz
spaces $\Lambda_\a^{p,\,q}$ by (1) and (2).\\

\noindent (ii). Concerning the Gleason solvability of the unit
ball for holomorphic mean Lipschitz spaces, discussed by Cho and
Zhu in ~\cite{7} (see, also, ~\cite{8}), we first remark that for
$0\,<\,\alpha_1\,<\,\alpha_2\,\leq\,1$,
$1\,\leq\,p,\,q_1\,<\,\infty$ and $q_2\,=\,\infty$, we have
\begin{equation} \label{Prop 3.6} \Lambda_{\alpha_2}^{p,\,q_1}\,\subset\,\Lambda_{\alpha_2}^p\,\subset\,\Lambda_{\alpha_1}^{p,\,q_1}\,\subset\,\Lambda_{\alpha_1}^p,
\end{equation}
by (i) and (ii) of ~\cite[Proposition 3.6]{7} (and all these
inclusions are continuous). Also, recall the roles of the indices
$\alpha,\, p$ and $q$ in the holomorphic mean Lipschitz space
$\Lambda_\alpha^{p,\,q}$ as follows. The index $p$ indicates the
basic $H^p$ norm that is used, the index $\a$ gives the order of
smoothness involved, and the index $q$ represents a rather subtle
correction to the order of smoothness; Proposition 3.6 of
~\cite{7} makes these remarks somewhat more precise. In fact, by
~\cite[Theorem 3.12]{7}, all these spaces are isomorphic. Further,
it is easy to see that ~\cite[Theorem 3.12]{7} also holds true for
the case $\alpha\,>\,0$, $s\,=\,q\,=\,\infty$ and $\beta\,=\,1/p$
with $1\,<\,p\,<\,\infty$; i.e.,
$\Lambda_\alpha\,\cong\,\Lambda_{\alpha - 1/p}$ (see, also,
~\cite[Theorem 4]{6}, used above). In ~\cite[Theorem 6]{6}, Chen
discussed certain inclusion relations between the classical
holomorphic Lipschitz spaces and holomorphic, integrated (i.e.,
$p\,<\,\infty$) Lipschitz spaces (see (1) above). The situation
for the classical holomorphic mean (integrated) Lipschitz spaces
of the unit ball will be similar due to the norm-equivalency
through the operators $J^\a$ and $\mathcal{R}^\a$ on these spaces,
and thus, we have
\begin{equation} \label{Thm 6_CZ} \Lambda_\alpha\,\subset\,\Lambda_\alpha^p\,\subset\,\Lambda_{\alpha - 1/p}
\end{equation} (and these inclusions are continuous). So, the space $\Lambda_\alpha^p$ is squeezed
between two classical holomorphic mean Lipschitz spaces, and by
(3), the space $\Lambda_\alpha^{p,\,q_1}$ is squeezed between two
classical holomorphic mean integrated Lipschitz spaces
$\Lambda_\alpha^p$. Now, applying the dead-end space method and
the interpolation method as in (i) above, we see that the unit
ball is Gleason solvable for all holomorphic mean (integrated)
Lipschitz spaces $\Lambda_\a^{p,\,q}$. Not only this, but we also
have the Gleason solvability of the unit ball for all holomorphic
mean Besov-Lipschitz spaces due to ~\cite[Theorems 1.2 and
1.3]{21} as argued at the end of (i) (this includes the special
case of the diagonal holomorphic Besov spaces $B_p =
\Lambda_{\frac{n}{p}}^{p, p}$ with $p < \infty$ and $\frac{n}{p} <
1$).

Concerning the Gleason solvability for Hardy-Sobolev spaces,
discussed by Cho and Zhu in ~\cite{7}, we have two cases: (A)
$p\,\in\,(1, 2]$, and (B) $p\,\in\,[2, \infty)$. For the case (A),
by (i) of Proposition 3.6, and (c) and (b) of Theorem 5.2 of
~\cite{7}, we have $\Lambda_\alpha^{p, p}\,\subset\,
H_\alpha^p\,\subset\,\Lambda_\alpha^{p,
2}\,\subset\,\Lambda_\alpha^p$ (and these inclusions are
continuous). For the case (B), by (i) of Proposition 3.6, and (a)
and (d) of Theorem 5.2 of ~\cite{7}, we have $\Lambda_\alpha^{p,
2}\,\subset\, H_\alpha^p\,\subset\,\Lambda_\alpha^{p,
p}\,\subset\,\Lambda_\alpha^p$ (and these inclusions are
continuous). So, as discussed in the previous paragraph, the space
$H_\a^p$ is squeezed between two Gleason solvable holomorphic mean
(integrated) Lipschitz spaces $\Lambda_\alpha^{p,\,q}$ and
$\Lambda_\alpha^{p}$. Now, applying the dead-end space method and
the interpolation method as in (i) above, we see that the unit
ball is Gleason solvable for all Hardy-Sobolev spaces. From the
definition of Besov-Hardy-Sobolev spaces, it is clear that the
unit ball is Gleason solvable for these spaces, since the unit
ball is Gleason solvable for Hardy-Sobolev spaces, and is Gleason
solvable for holomorphic Besov-Lipschitz spaces by (i).
Alternatively, we have the following analogue of (3) from the
definition of Besov-Hardy-Sobolev spaces
\begin{equation} \label{Prop 3.6_BHS} H_{\alpha_2}^{p,\,q_1}\,\subset\,H_{\alpha_2}^p\,\subset\,H_{\alpha_1}^{p,\,q_1}\,\subset\,H_{\alpha_1}^p
\end{equation}
(and these inclusions are continuous), showing that the space
$H_\a^{p,\,q_1}$ is squeezed between two Hardy-Sobolev spaces. It
is easy to obtain an analogue of ~\cite[Theorem 3.13]{7}, to show
that all these spaces are isomorphic. Now, applying the dead-end
space method the interpolation method as in (i) above, we see that
the unit ball is Gleason solvable for all
Besov-Hardy-Sobolev spaces since so are Hardy-Sobolev spaces.\\

\noindent (iii). We recall that a special class of mixed-norm
spaces from ~\cite{7} is actually the class of mixed-norm Bergman
spaces, namely, Besov-Lipschitz spaces $\Lambda_\a^{p,\,q}$,
discussed in ~\cite{21}. So, by (i), the unit ball is Gleason
solvable for such mixed-norm spaces. Moreover, weighted Bergman
spaces are closely related with Hardy-Sobolev spaces as well as
with Besov-Lipschitz spaces through the radial fractional
differential operators with suitable exponents. Hence the unit
ball is also Gleason solvable for weighted Bergman spaces since it
is Gleason solvable for Hardy-Sobolev spaces by (ii) and for
Besov-Lipschitz spaces by (i).\\

\noindent (iv). Next, we obtain the Gleason solvability for some
special classes of holomorphic (mean) Besov-Lipschitz spaces,
discussed by Jevti\'{c} and Pavlovi\'{c} in ~\cite{21}. For this,
we note that for $0\,<\,\alpha\,<\,n$ and $\phi$ an increasing,
positive, continuous function on $(0, \,1]$, the space
$\Lambda_{n, \phi}^p$ contains the Hardy-Sobolev space $H_n^p$,
and the integrated space $\Lambda_{n, \phi}^{p, q}$ contains the
Besov-Hardy-Sobolev space $H_n^{p, q}$. Further, if the function
$x\,\mapsto \phi(x)/x^\alpha$ increases on $(0, \,1]$, then
$\Lambda_{n, \phi}^p\,\subset\,\Lambda_\alpha^p$; in fact, we have
$\Lambda_{n, \phi}^p\,\cong\,\textrm{Lip}_{n, \phi}^p$, by
Corollary 1.9, (1.12) and (1.16) of ~\cite{21}. The same
situations hold true for the integrated spaces $\Lambda_{n,
\phi}^{p, q}$, $\Lambda_\alpha^{p, q}$ and $\textrm{Lip}_{n,
\phi}^{p, q}$, by Corollary 1.10, (1.11), (1.18) and (1.19) of
~\cite{21}. The space $\Lambda_{n, \phi}^p$ (resp., $\Lambda_{n,
\phi}^{p, q}$) is squeezed between the Hardy-Sobolev space $H_n^p$
(resp., the Besov-Hardy-Sobolev space $H_n^{p,\,q}$) and the
holomorphic mean Lipschitz space $\Lambda_\a^p$ (resp., the
holomorphic mean Besov-Lipschitz space $\Lambda_\a^{p,\,q}$).
Since the (Besov-)Hardy-Sobolev spaces and holomorphic mean
(Besov-)Lipschitz spaces are Gleason solvable, we see that the
unit ball is Gleason solvable for these special classes of
holomorphic (mean) Besov-Lipschitz spaces.\\

\noindent {\it Proof of Theorem 3.1 (B).} Now, using the dead-end
space method and the interpolation method, we establish the corona
theorem for the same Banach spaces of holomorphic functions on
bounded symmetric domains in $\cc^n$ (including the unit ball and
the polydisc). We do not repeat the central arguments for these
spaces, discussed in the second preceding paragraph to the proof
of Theorem 3.1 (A). However we remark that if $F$ is the middle
space in the inclusion relation for which we wish to establish the
corona theorem for a given corona data in that space, then it is
clear that we have the same corona data in $F_0$ (and hence, we
have a corona data in $F_\infty$, as we have $F_0 \cong
F_\infty$). Now, since the corona theorem holds true in both
$F_\infty$ and $F_0$ (also, in $E_\infty$, $E_r$ and $E_0$ as
well) by Theorem 2.2 (due to
Theorem 2.1), it also holds true in $F$ as well. $\hfill \Box$\\

\noindent {\bf Open Problems.} Since the main source of main
results is author's previous paper on affirmative solutions of
these two problems for holomorphic Lipschitz algebras on
``certain" bounded pseudoconvex/poly domains in $\cc^n$, it is
interesting to extend Theorem 3.1 to other suitable domains $D$
(e.g., ``polydomains", considered by author in ~\cite{28}, and
more specifically, a product of Reinhardt domains) in $\cc^n$. In
particular, it is of interest to know whether Theorem 3.1 (ii) to
(iv) hold true for the polydisc in $\cc^n$.\\

\noindent {\bf Funding.} The author did not receive financial
support from any organization for the submitted work.\\

\noindent {\bf Conflict of interest.} The author declares that he
has no conflicts of interest.

\end{document}